\documentclass{conm-p-l}

\newtheorem{theorem}{Theorem}[section]

\theoremstyle{definition}
\newtheorem{definition}[theorem]{Definition}

\theoremstyle{remark}

\numberwithin{equation}{section}

\newcommand{\LL}{{\mathcal L}}
\newcommand{\tla}{\tilde{\lambda}}

\newcommand{\hk}{\hat{k}}
\newcommand{\Z}{\mathbb{Z}^2/\{0\}}

\newcommand{\e}{\epsilon}

\renewcommand{\k}{\kappa}

\newcommand{\Ga}{\Gamma}

\newcommand{\Dl}{\Delta}

\newcommand{\Sg}{\Sigma}

\newcommand{\pa}{\partial}

\newcommand{\la}{\lambda}

\newcommand{\nid}{\noindent}

\newcommand{\om}{\omega}
\newcommand{\Om}{\Omega}

\newcommand{\lag}{\langle}
\newcommand{\rag}{\rangle}

\begin{document}

\title{On 2D Euler Equations: III. A Line Model}

\author{Yanguang (Charles)  Li}
\address{Department of Mathematics, University of Missouri, 
Columbia, MO 65211}
\curraddr{}
\email{cli@math.missouri.edu}
\thanks{}


\subjclass{Primary 76, 35; Secondary 34, 37.}
\date{}



\begin{abstract}
The spectral theorem of the linear 2D Euler operator in Sobolev spaces 
is presented as a corollary of the spectral theorem in $\ell_2$ space 
in \cite{Li00}. Study on the (dashed) line model introduced in 
\cite{Li01} is continued. Specifically, invariant manifolds for the 
line model is established. The corresponding line model for 2D Navier-Stokes 
equation is also introduced.
\end{abstract}

\maketitle








\section{Introduction}

To understand the nature of turbulence, we select 2D Euler equation under 
periodic boundary condition as our primary example to study. 2D 
Navier-Stokes equation at high Reynolds number is regarded as a singularly 
perturbed 2D Euler equation. That is, we are interested in studying the 
zero viscosity limit problem.

To begin an infinite dimensional dynamical system study, we consider a 
simple fixed point and study the spectrum of the linear 2D Euler operator 
in \cite{Li00}. The spectral theorem in $\ell_2$ space is proved. As a 
corollary of the spectral theorem in $\ell_2$ space, we will present the 
spectral theorem in Sobolev spaces in this article. Sobolev spaces are of 
more interest to us, since we are interested in understanding the invariant 
manifolds of 2D Euler equation at the fixed point. The main obstacle toward 
proving the invariant manifold theorem is that the nonlinear term is 
non-Lipschitzian. In \cite{Li01}, a (dashed) line model is introduced 
to understand the invariant manifold structure of 2D Euler equation. At 
a special parameter value, the explicit expression of the invariant manifolds 
of the dashed line model can be calculated. The stable and unstable manifolds 
are two dimensional ellipsoidal surfaces, and together they form a lip-shape 
hyperbolic structure. Such structure appears to be robust with repsect to 
the parameter. In this article, we will prove the existence of invariant 
manifolds for the line model.

Another more exciting development is the discovery of a Lax pair for 2D 
Euler equation \cite{Li01a}. From the Lax pair, we have obtained a Darboux 
transformation for the 2D Euler equation \cite{LY00}. In principle, explicit 
expressions of the hyperbolic structures can be obtained from Darboux 
transformations \cite{Li00a}.

Hyperbolic structures are the source of chaos when the system is under 
perturbations. The corresponding line model for 2D Navier-Stokes equation 
with temporally periodic forcing is a singular perturbation of the line 
model for 2D Euler equation. Numerical simulations on these line models 
are conducted. 

The article is organized as follows: In section 2, we present the 
spectral theorem in Sobolev spaces. In section 3, we present the  
invariant manifold result on the line model. In section 4, we 
introduce the corresponding line model for 2D Navier-Stokes equation.

\section{The Spectral Theorem of the Linear 2D Euler Operator in 
Sobolev Spaces}

Consider the 2D Euler equation in vorticity form
\begin{equation}
{\pa \Om \over \pa t} + \{ \Psi, \Om \} = 0 \ ,
\label{euler}
\end{equation}
where $\Om$ is the vorticity, $\Psi$ is the stream function, 
$\Om = \Dl \Psi$, $\Dl$ is the 2D Laplacian, and
the bracket $\{\ ,\ \}$ is defined as
\[
\{ f, g\} = (\pa_x f) (\pa_y g) - (\pa_y f) (\pa_x g) \ .
\]
Expanding $\Om$ into Fourier series,
\[
\Om =\sum_{k\in \Z} \om_k \ e^{ik\cdot X}\ ,
\]
where $\om_{-k}=\overline{\om_k}\ $, $k=(k_1,k_2)^T$, 
and $X=(x,y)^T$. The 2D Euler equation
can be rewritten as 
\begin{equation}
\dot{\om}_k = \sum_{k=p+q} A(p,q) \ \om_p \om_q \ ,
\label{Keuler}
\end{equation}
where $A(p,q)$ is given by,
\begin{eqnarray}
A(p,q)&=& {1\over 2}[|q|^{-2}-|p|^{-2}](p_1 q_2 -p_2 q_1) \nonumber \\
\label{Af} \\      
      &=& {1\over 2}[|q|^{-2}-|p|^{-2}]\left | \begin{array}{lr} 
p_1 & q_1 \\ p_2 & q_2 \\ \end{array} \right | \ , \nonumber
\end{eqnarray}
where $|q|^2 =q_1^2 +q_2^2$ for $q=(q_1,q_2)^T$, similarly for $p$.
Denote $\{ \om_k \}_{k\in \Z}$ by $\om$. We consider the simple fixed point 
$\om^*$:
\begin{equation}
\om^*_p = \Ga,\ \ \ \om^*_k = 0 ,\ \mbox{if} \ k \neq p \ \mbox{or}\ -p,
\label{fixpt}
\end{equation}
of the 2D Euler equation (\ref{Keuler}), where 
$\Ga$ is an arbitrary complex constant. 
The {\em{linearized two-dimensional Euler equation}} at $\om^*$ is given by,
\begin{equation}
\dot{\om}_k = A(p,k-p)\ \Ga \ \om_{k-p} + A(-p,k+p)\ \bar{\Ga}\ \om_{k+p}\ .
\label{LE}
\end{equation}
\begin{definition}[Classes]
For any $\hk \in \Z$, we define the class $\Sg_{\hk}$ to be the subset of 
$\Z$:
\[
\Sg_{\hk} = \bigg \{ \hk + n p \in \Z \ \bigg | \ n \in \mathbb{Z}, \ 
\ p \ \mbox{is specified in (\ref{fixpt})} \bigg \}.
\]
\label{classify}
\end{definition}
\nid
See Fig.\ref{class} for an illustration of the classes. 
According to the classification 
defined in Definition \ref{classify}, the linearized two-dimensional Euler 
equation (\ref{LE}) decouples into infinitely many {\em{invariant subsystems}}:
\begin{eqnarray}
\dot{\omega}_{\hat{k} + np} &=& A(p, \hat{k} + (n-1) p) 
     \ \Gamma \ \omega_{\hat{k} + (n-1) p} \nonumber \\
& & + \ A(-p, \hat{k} + (n+1)p)\ 
     \bar{\Gamma} \ \omega_{\hat{k} +(n+1)p}\ . \label{CLE}
\end{eqnarray}
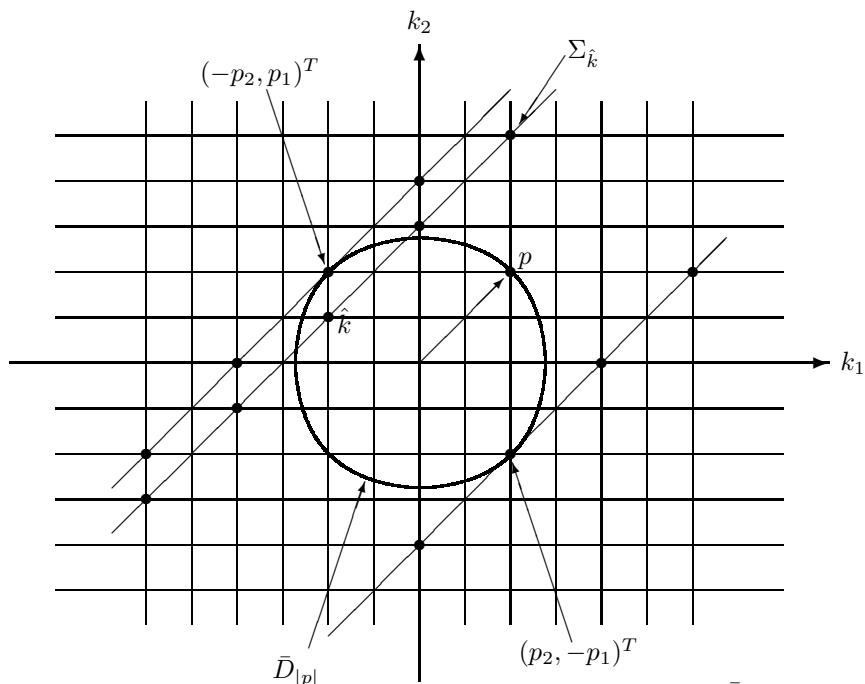
\begin{figure}[ht]
  \begin{center}
    \leavevmode
      \setlength{\unitlength}{2ex}
  \begin{picture}(36,27.8)(-18,-12)
    \thinlines
\multiput(-12,-11.5)(2,0){13}{\line(0,1){23}}
\multiput(-16,-10)(0,2){11}{\line(1,0){32}}
    \thicklines
\put(0,-14){\vector(0,1){28}}
\put(-18,0){\vector(1,0){36}}
\put(0,15){\makebox(0,0){$k_2$}}
\put(18.5,0){\makebox(0,0)[l]{$k_1$}}
\qbezier(-5.5,0)(-5.275,5.275)(0,5.5)
\qbezier(0,5.5)(5.275,5.275)(5.5,0)
\qbezier(5.5,0)(5.275,-5.275)(0,-5.5)
\qbezier(0,-5.5)(-5.275,-5.275)(-5.5,0)
    \thinlines
\put(4,4){\circle*{0.5}}
\put(0,0){\vector(1,1){3.7}}
\put(4.35,4.35){$p$}
\put(4,-4){\circle*{0.5}}
\put(8,0){\circle*{0.5}}
\put(-8,0){\circle*{0.5}}
\put(-8,-2){\circle*{0.5}}
\put(-12,-4){\circle*{0.5}}
\put(-12,-6){\circle*{0.5}}
\put(-4,2){\circle*{0.5}}
\put(-4,4){\circle*{0.5}}
\put(0,6){\circle*{0.5}}
\put(0,8){\circle*{0.5}}
\put(4,10){\circle*{0.5}}
\put(12,4){\circle*{0.5}}
\put(0,-8){\circle*{0.5}}
\put(-4,-12){\line(1,1){17.5}}
\put(-13.5,-7.5){\line(1,1){19.5}}
\put(-13.5,-5.5){\line(1,1){17.5}}
\put(-3.6,1.3){$\hat{k}$}
\put(-7,12.1){\makebox(0,0)[b]{$(-p_2, p_1)^T$}}
\put(-6.7,12){\vector(1,-3){2.55}}
\put(6.5,13.6){\makebox(0,0)[l]{$\Sg_{\hat{k}}$}}
\put(6.4,13.5){\vector(-2,-3){2.0}}
\put(7,-12.1){\makebox(0,0)[t]{$(p_2, -p_1)^T$}}
\put(6.7,-12.25){\vector(-1,3){2.62}}
\put(-4.4,-13.6){\makebox(0,0)[r]{$\bar{D}_{|p|}$}}
\put(-4.85,-12.55){\vector(1,3){2.45}}
\end{picture}
\end{center}
\caption{An illustration of the classes $\Sg_{\hk}$ and the disk 
$\bar{D}_{|p|}$.}
\label{class}
\end{figure}
Let $\LL_{\hk}$ be the linear operator defined by the right hand side of 
(\ref{CLE}), and $H^s$ be the Sobolev space where $s \geq 0$ is an integer and 
$H^0=\ell_2$. 
\begin{theorem}
The eigenvalues of the linear operator $\LL_{\hk}$ in $H^s$ are of 
four types: real pairs ($c, -c$), purely imaginary pairs ($id, -id$), 
quadruples ($\pm c \pm id$), and zero eigenvalues.
\end{theorem}
Proof. The same proof as in \cite{Li00} works here. QED

\nid
The eigenvalues can be computed through continued fractions \cite{Li00}.
\begin{definition}[The Disk]
The disk of radius $| p |$ in $\Z$, denoted by
$\bar{D}_{| p |}$, is defined as
\[ 
 \bar{D}_{| p |} = \bigg \{ k \in \Z \ \bigg| 
     \ | k | \leq | p | \bigg \} \, .
\]
\end{definition}
\nid
See Fig.\ref{class} for an illustration. 
\begin{theorem}[The Spectral Theorem] We have the following claims on 
the spectrum of the linear operator $\LL_{\hk}$:
\begin{enumerate}
\item If $\Sg_{\hat{k}} \cap \bar{D}_{|p|} = \emptyset$, then the entire
$H^s$ spectrum of the linear operator $\LL_{\hk}$ 
is its continuous spectrum. See Figure \ref{splb}, where
$b= - \frac{1}{2}|\Gamma | |p|^{-2} 
\left|
  \begin{array}{cc}
p_1 & \hat{k}_1 \\
p_2 & \hat{k}_2
  \end{array}
\right| \ .$
That is, both the residual and the point spectra of $\LL_{\hk}$ are empty.
\item If $\Sg_{\hat{k}} \cap \bar{D}_{|p|} \neq \emptyset$, then the entire
essential $H^s$ spectrum of the linear operator $\LL_{\hk}$ is its 
continuous spectrum. 
That is, the residual 
spectrum of $\LL_{\hk}$ is empty. The point 
spectrum of $\LL_{\hk}$ is symmetric with respect to both real and 
imaginary axes. 
See Figure \ref{spla2}.
\end{enumerate}
\label{SST}
\end{theorem}
Proof. The same proof as in \cite{Li00} can be carried through here for 
$H^s$, with the following minor modifications:
\begin{enumerate}
\item In the proof of Theorem VI.1 on page 747 of \cite{Li00}, simply 
replace the $\ell_2$-norm by $H^s$-norm.
\item In the proof of Theorem VI.3 on page 750 of \cite{Li00}, the inner
product $\lag \ , \ \rag$ should still be an $\ell_2$ inner product.
\item In the proof of Theorem VI.4 on page 751 of \cite{Li00}, $\ell_1$, 
$\ell_2$, and  $\ell_\infty$ should be replaced by the Sobolev spaces 
$W^{s,1}$, $W^{s,2}$ ($=H^s$ in our notation), and $W^{s,\infty }$. In the 
expression (VI.60) of $f_n$ on page 753, one has 
\begin{eqnarray*}
 n^sf_n &=& \frac{2 \left[ 1-w_*^{-4} \right]}{W_0} 
       \ \sum^{n-1}_{j=0} \left[ (n-j)^s w_*^{n-j} + (n-j)^s 
(-w_*)^{n-j} \right] \\
& & \bigg ( \frac{n^s}{(n-j)^s (j+2)^s} \bigg )(j+2)^sy_{j+2} 
     - \frac{2 \left[ 1-w_*^4 \right]}{W_0} 
       \sum^{\infty}_{j=n} \\
& & \left[ w_*^{j-n} + (-w_*)^{j-n} \right]
\bigg ( \frac{n^s}{(j+2)^s} \bigg )(j+2)^s y_{j+2} \, ,  
\end{eqnarray*}
where both $\bigg ( \frac{n^s}{(n-j)^s (j+2)^s} \bigg )$ and 
$\bigg ( \frac{n^s}{(j+2)^s} \bigg )$ are bounded in $n$ and $j$, and 
the Riesz convexity theorem can be applied.
\end{enumerate}
The proof of the theorem is completed. QED

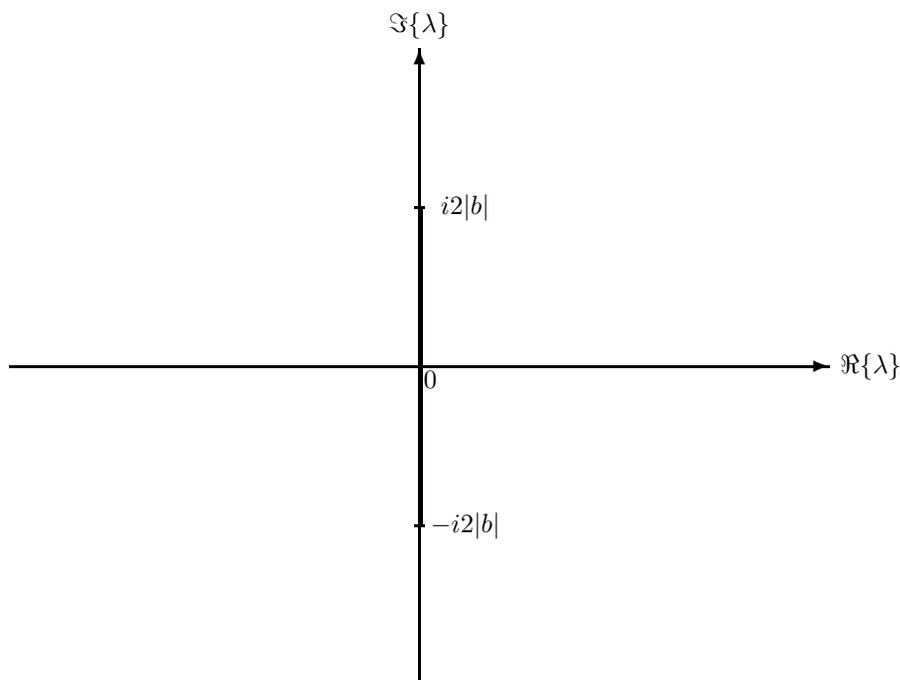
\begin{figure}[ht]
  \begin{center}
    \leavevmode
      \setlength{\unitlength}{2ex}
  \begin{picture}(36,27.8)(-18,-12)
    \thicklines
\put(0,-14){\vector(0,1){28}}
\put(-18,0){\vector(1,0){36}}
\put(0,15){\makebox(0,0){$\Im \{ \la \}$}}
\put(18.5,0){\makebox(0,0)[l]{$\Re \{ \la \}$}}
\put(0.1,-7){\line(0,1){14}}
\put(.2,-.2){\makebox(0,0)[tl]{$0$}}
\put(-0.2,-7){\line(1,0){0.4}}
\put(-0.2,7){\line(1,0){0.4}}
\put(2.0,-6.4){\makebox(0,0)[t]{$-i2|b|$}}
\put(2.0,7.6){\makebox(0,0)[t]{$i2|b|$}}
\end{picture}
  \end{center}
\caption{The spectrum of $\LL_{\hk}$ in case (1).}
\label{splb}
\end{figure}
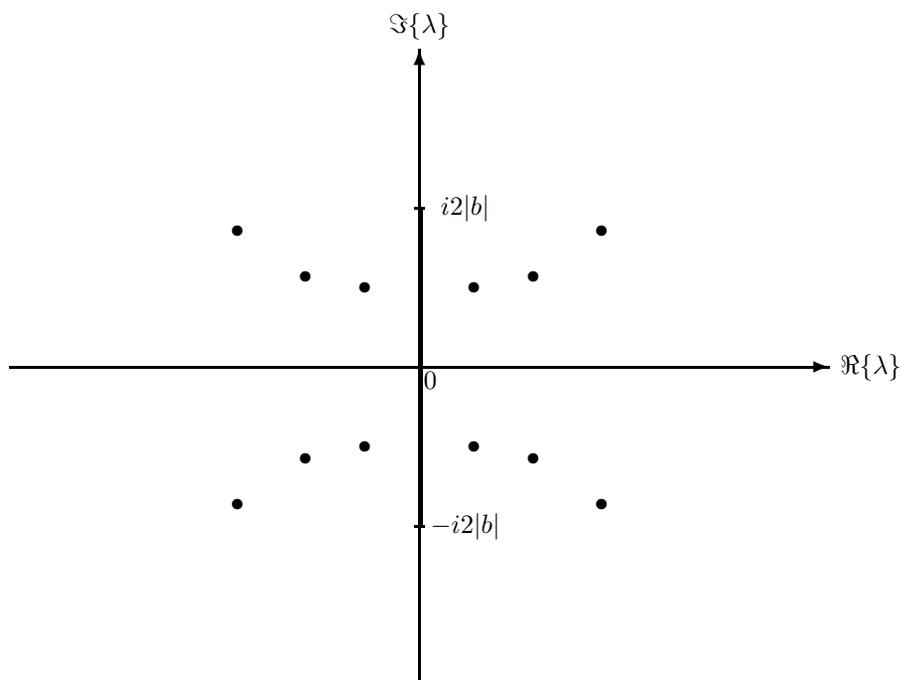
\begin{figure}[ht]
  \begin{center}
    \leavevmode
      \setlength{\unitlength}{2ex}
  \begin{picture}(36,27.8)(-18,-12)
    \thicklines
\put(0,-14){\vector(0,1){28}}
\put(-18,0){\vector(1,0){36}}
\put(0,15){\makebox(0,0){$\Im \{ \la \}$}}
\put(18.5,0){\makebox(0,0)[l]{$\Re \{ \la \}$}}
\put(0.1,-7){\line(0,1){14}}
\put(.2,-.2){\makebox(0,0)[tl]{$0$}}
\put(-0.2,-7){\line(1,0){0.4}}
\put(-0.2,7){\line(1,0){0.4}}
\put(2.0,-6.4){\makebox(0,0)[t]{$-i2|b|$}}
\put(2.0,7.6){\makebox(0,0)[t]{$i2|b|$}}
\put(2.4,3.5){\circle*{0.5}}
\put(-2.4,3.5){\circle*{0.5}}
\put(2.4,-3.5){\circle*{0.5}}
\put(-2.4,-3.5){\circle*{0.5}}
\put(5,4){\circle*{0.5}}
\put(-5,4){\circle*{0.5}}
\put(5,-4){\circle*{0.5}}
\put(-5,-4){\circle*{0.5}}
\put(8,6){\circle*{0.5}}
\put(-8,6){\circle*{0.5}}
\put(8,-6){\circle*{0.5}}
\put(-8,-6){\circle*{0.5}}
\end{picture}
  \end{center}
\caption{The spectrum of $\LL_{\hk}$ in case (2).}
\label{spla2}
\end{figure}
\begin{figure}[ht]
  \begin{center}
    \leavevmode
      \setlength{\unitlength}{2ex}
  \begin{picture}(36,27.8)(-18,-12)
    \thicklines
\put(0,-14){\vector(0,1){28}}
\put(-18,0){\vector(1,0){36}}
\put(0,15){\makebox(0,0){$\Im \{ \la \}$}}
\put(18.5,0){\makebox(0,0)[l]{$\Re \{ \la \}$}}
\put(2.4,3.5){\circle*{0.5}}
\put(-2.4,3.5){\circle*{0.5}}
\put(2.4,-3.5){\circle*{0.5}}
\put(-2.4,-3.5){\circle*{0.5}}  
\put(0.1,-10){\line(0,1){20}}
\put(.2,-.2){\makebox(0,0)[tl]{$0$}}
\put(-0.2,-10){\line(1,0){0.4}}
\put(-0.2,10){\line(1,0){0.4}}
\put(2.0,-9.4){\makebox(0,0)[t]{$-i2|b|$}}
\put(2.0,10.6){\makebox(0,0)[t]{$i2|b|$}}
\end{picture}
  \end{center}
\caption{The spectrum of $\LL_{\hk}$ with $\hk = (-3,-2)^T$, when $p=(1,1)^T$.}
\label{figev}
\end{figure}
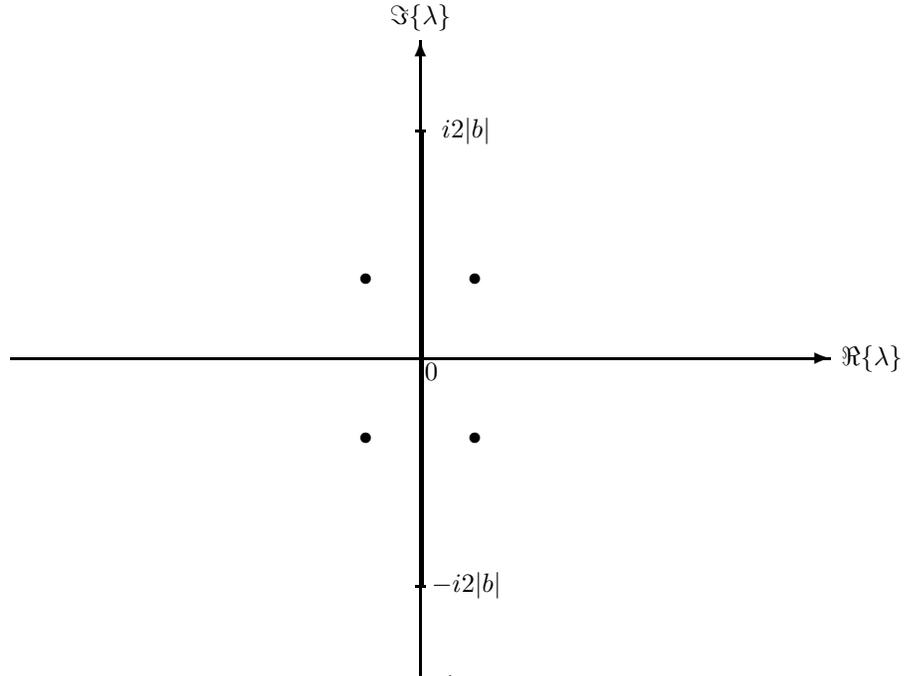
\begin{figure}[ht]
  \begin{center}
    \leavevmode
      \setlength{\unitlength}{2ex}
  \begin{picture}(36,27.8)(-18,-12)
    \thinlines
\multiput(-12,-11.5)(2,0){13}{\line(0,1){23}}
\multiput(-16,-10)(0,2){11}{\line(1,0){32}}
    \thicklines
\put(0,-14){\vector(0,1){28}}
\put(-18,0){\vector(1,0){36}}
\put(0,15){\makebox(0,0){$k_2$}}
\put(18.5,0){\makebox(0,0)[l]{$k_1$}}
\qbezier(-2.75,0)(-2.6375,2.6375)(0,2.75)
\qbezier(0,2.75)(2.6375,2.6375)(2.75,0)
\qbezier(2.75,0)(2.6375,-2.6375)(0,-2.75)
\qbezier(0,-2.75)(-2.6375,-2.6375)(-2.75,0)
    \thinlines
\put(2,2){\circle*{0.5}}
\put(0,0){\vector(1,1){1.85}}
\put(2.275,2.275){$p$}
\put(-12,-10){\circle*{0.5}}
\put(-10,-8){\circle*{0.5}}
\put(-8,-6){\circle*{0.5}}
\put(-6,-4){\circle*{0.5}}
\put(-4,-2){\circle*{0.5}}
\put(-2,0){\circle*{0.5}}
\put(0,2){\circle*{0.5}}
\put(2,4){\circle*{0.5}}
\put(4,6){\circle*{0.5}}
\put(6,8){\circle*{0.5}}
\put(8,10){\circle*{0.5}}
\put(-14,-12){\line(1,1){24}}
\put(-5.6,-5.4){$\hat{k}$}
\put(-4.4,-13.6){\makebox(0,0)[r]{$\bar{D}_{|p|}$}}
\put(-4.85,-12.55){\vector(1,3){3.4}}
\end{picture}
\end{center}
\caption{The collocation of the modes in the line model.}
\label{model}
\end{figure}
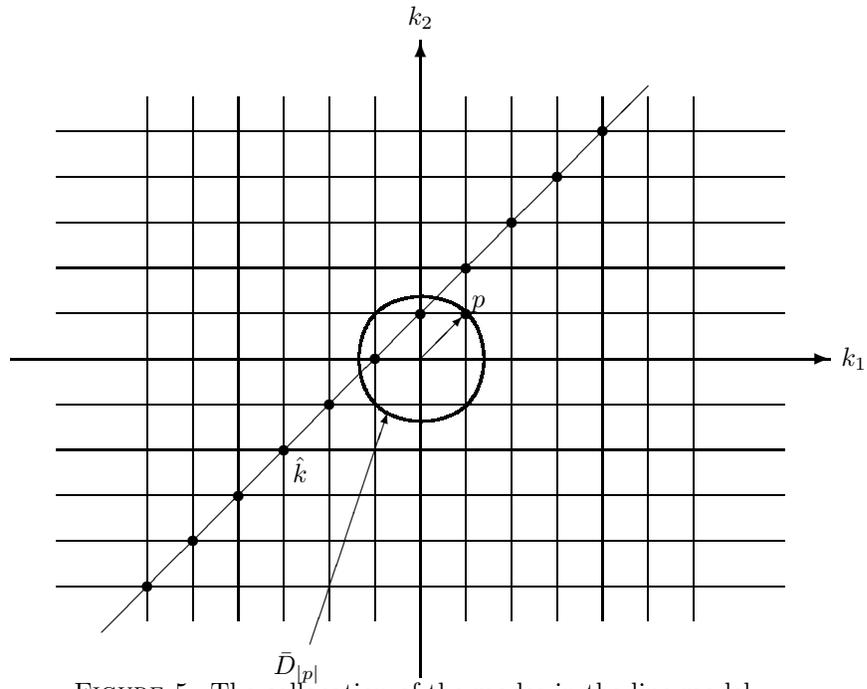

\section{Invariant Manifolds of the Line Model}

To simplify our study, we study only the case when $\om_k$ is real, $\forall 
k \in \Z$, i.e. we only study the cosine transform of the vorticity, 
\[
\Om = \sum_{k \in \Z} \om_k \cos (k \cdot X)\ ,
\]
and the 2D Euler equation (\ref{euler}) preserves the cosine 
transform. To further simplify our study, we will study a concrete 
line model based upon the fixed point (\ref{fixpt}) with the 
mode $p=(1,1)^T$ and parametrized by $\Ga$.
When $\Ga \neq 0$, the fixed point has $4$ eigenvalues which form a 
quadruple. These four eigenvalues appear in the invariant 
linear subsystem labeled by $\hk = (-3,-2)^T$. We computed the eigenvalues 
through continued fractions, one of them is \cite{Li00}:
\begin{equation}
\tla=2 \lambda / | \Gamma | = 0.24822302478255 \ + \ i \ 0.35172076526520\ .
\label{evun}
\end{equation}
See Figure \ref{figev} for an illustration. The essential spectrum 
(= continuous spectrum) of $\LL_{\hk}$ with $\hk = (-3,-2)^T$ is the segment 
on the imaginary axis shown in Figure \ref{figev}, where 
$b = -\frac{1}{4} \Ga$. The essential spectrum 
(= continuous spectrum) of the linear 2D Euler operator at this fixed point 
is the entire imaginary axis. The line model is a Galerkin truncation with 
the modes on the line $\{ \hk + n p , \ n \in \mathbb{Z} \}$ and $p$, where 
$\hk = (-3,-2)^T$ and $p=(1,1)^T$. See Figure \ref{model} for an 
illustration of the modes in this model, which has the 
line nature leading to the name of the model. The line model is designed to 
model the hyperbolic structure in the neighborhood of the fixed point. For 
simplicity of presentation, we use the abbreviated notations,
\[
  \omega_n = \omega_{\hat{k}+np} \, , \ \  
  A_n = A(p,\hat{k}+np) \, , \ \  
  A_{m,n} = A(\hat{k}+mp,\hat{k}+np) \, .  
\]
The {\em{line model}} is,
\begin{eqnarray}
\dot{\omega}_n &=& A_{n-1} \omega_p
\omega_{n-1} - A_{n+1} \omega_p \omega_{n+1} \ , 
\nonumber \\ 
\label{rdlm} \\
\dot{\omega}_p &=& - \sum_{n \in Z}
A_{n-1,n} \omega_{n-1} \omega_n \, . \nonumber
\end{eqnarray}
We also use $\om^*$ to denote the fixed point of the line model 
$\{ \om_p = \Ga; \ \om_n = 0, \ \forall n \in \mathbb{Z} \}$. The linearized 
line model at the fixed point $\om^*$ is the same with the invariant 
subsystem (\ref{CLE}) with $\hk = (-3,-2)^T$ and $p=(1,1)^T$. In the 
neighborhood of the fixed point $\om^*$, the line model can be rewritten as
\begin{equation}
\dot{\om} = L\om +Q(\om)\ , 
\label{mlm}
\end{equation}
where
\begin{eqnarray*}
& & \om = (\om_p, \quad \om_n \quad (n \in \mathbb{Z}))\ , \\
& & [L\om]_n = A_{n-1} \Ga \om_{n-1} - A_{n+1} \Ga \om_{n+1} \ , \\
& & [L\om]_p =0 \ , \\
& & [Q(\om)]_n = A_{n-1} \om_p \om_{n-1} - A_{n+1} \om_p \om_{n+1} \ , \\
& & [Q(\om)]_p = - \sum_{n \in \mathbb{Z}}A_{n-1,n}\om_{n-1} \om_{n}\ . 
\end{eqnarray*}
The spectrum of $L$ is given in Figure \ref{figev}. We have
\begin{theorem}[Invariant Manifold Theorem]
The line of fixed points $\om^*$ parametrized by $\Ga$ of the line model has 
codimension 2 smooth center-unstable and center-stable manifolds, and 
codimension 4 smooth center manifold, in $H^s$ ($s \geq 1$).
\end{theorem}
Proof. The spectrum of $L$ in $H^s$ ($s \geq 1$) is given in Figure 
\ref{figev}. $Q(\om)$ is quadratic in $\om$ for $H^s$ ($s \geq 1$) is 
a Banach algebra. In a small neighborhood of $\om^*$, the existence of 
smooth center-unstable, center-stable, and center manifolds follows from 
standard arguments. QED

\section{A Line Model for 2D Navier-Stokes Equation}

We consider the 2D Navier-Stokes equation with temporally periodic forcing,
as a singular perturbation of the 2D Euler equation (\ref{euler}),
\begin{equation}
{\pa \Om \over \pa t} + \{ \Psi, \Om \} = \e [\Dl \Om + f(t,x,y)] \ ,
\label{nse}
\end{equation}
where $\e = 1/\mbox{Re}$ is the inverse of Reynolds number, and $f(t,x,y)$ 
is periodic in $t$, periodic in $x$ and $y$ of period $2\pi$, and of 
spatial mean $0$. The corresponding line model is
\begin{eqnarray}
\dot{\omega}_n &=& A_{n-1} \omega_p
\omega_{n-1} - A_{n+1} \omega_p \omega_{n+1} +
\e [-\k_n^2 \om_n +f_n(t)]\ , \nonumber \\ 
\label{nslm} \\
\dot{\omega}_p &=& - \sum_{n \in Z}
A_{n-1,n} \omega_{n-1} \omega_n +
\e [-\k_p^2 \om_p +f_p(t)]\, , \nonumber
\end{eqnarray}
where $f_n$ and $f_p$ are periodic functions of $t$, and 
\[
\k_n = | \hk +np|\ , \quad \k_p =|p|\ , \quad \hk = (-3,-2)^T\ ,
\quad p=(1,1)^T \ .
\]

\end{document}